\documentclass[14pt]{extarticle}
\usepackage{standalone}
\usepackage{amsmath}
\usepackage{amssymb}
\usepackage{amsbsy}
\usepackage{amsfonts}
\usepackage{eucal}
\usepackage{titlesec}
\usepackage{caption}
\usepackage{tikz-cd}
\usepackage[all,cmtip]{xy}
\usepackage{concrete}
\usepackage{amsthm}
\usepackage{mathrsfs}
\usepackage{geometry,graphicx,color}
\usepackage{mathtext}
\usepackage{hyperref}
\usepackage{booktabs}
\usepackage{fancybox}
\theoremstyle{plain}
\newtheorem{theorem}{Theorem}
\newtheorem{lemma}[theorem]{Lemma}

\newtheorem{definition}[theorem]{Definition}

\newcommand{\xar}[1]{\ensuremath{\xrightarrow{#1}}}
\newcommand{\mf}[1]{\pmb #1 }
\newcommand{\mc}[1]{\mathcal{#1}}

\newcommand{\R}{\mathbb{R}}

\newcommand{\Z}{\mathbb{Z}}
\newcommand{\C}{\mathbb{C}}

\newcommand{\on}[1]{\operatorname{#1}}
\newcommand{\la}{\langle}
\newcommand{\ra}{\rangle}	
\newcommand{\ol}[1]{\overline{#1}}
\newcommand{\ul}[1]{\underline{#1}}
\newcommand{\op}{\mathrm{op}}
\newcommand{\N}{{\pmb  \Delta}} 
\newcommand{\ca}{\circlearrowright}

\usepackage{chngcntr}
\titleformat{\subsubsection}[runin]{\normalsize\bfseries}{\textbf{(\arabic{section}.\arabic{subsection}.\arabic{subsubsection}})}{0.5em}{}[]
\newcommand{\p}[1]{\subsection{#1}}
\titleformat{\subsection}[runin]{\normalsize\bfseries}{\textbf{(\arabic{section}.\arabic{subsection}})}{0.5em}{}[]
\titleformat{\section}[block]{\normalsize\bfseries}{\textbf{\arabic{section}.}}{0.5em}{}[]

\title{$K(\Z,2)$ out of circular permutations}
\author{Nikolai Mn\"ev\thanks{PDMI RAS;  Chebyshev Laboratory, SPbSU. Research is supported by the Russian Science Foundation grant 19-71-30002.} \\  \href{mailto:mnev@pdmi.ras.ru}{mnev@pdmi.ras.ru} }

\begin{document}
\maketitle

\begin{abstract}
	We  discuss $\pmb{SC}_*$, a simplicial homotopy model of
	$K(\Z,2)$ constructed from circular permutations. In any dimension,
	the number of simplices in the model is finite. The complex $\pmb{SC}_*$
	naturally manifests as a simplicial set representing ``minimally"
	triangulated circle bundles over simplicial bases. On the other hand,
	existence of the homotopy equivalence  $|\pmb{SC}_*| \approx B(U(1)) \approx K(\Z,2)$ appears to
	be a canonical fact from the foundations of the theory of crossed simplicial groups.
\end{abstract}
\tableofcontents
\section{Introduction}
This note essentially continues the discussion from \cite{MnevMin}. In
that note (\cite[$\S\S$ 3.6, 3.7]{MnevMin}), we identify circular
permutations of $n+1$ ordered elements with ``minimal" semi-simplicial
triangulations of trivial circle bundles over ordered base
$n$-simplices. Any semi-simplicial triangulation of a circle bundle is
non-canonically combinatorially concordant to a minimal triangulation
(i.e., having minimal triangulations over all the simplices of the
same base complex), and the simplicial set $\pmb{SC}_*$ of circular
permutations naturally represents minimally triangulated circle
bundles over semi-simplicial complexes. Such triangulations
functorially (via Kan's second derived subdivision
$\operatorname{Sd}_2$) have the structure of a classical simplicial PL
triangulation. However, the \textit{minimal} triangulations exist only
in the semi-simplicial category. The value (if it exists) of the above
constructions lies in their very discrete form of the Weil-Kostant
correspondence for triangulated circle bundles (\cite[Theorem
1]{MnevMin}). Namely, a circle bundle over a given simplicial complex
$B$ can be (semi-simplicially) triangulated with base $B$ if and only
if its Chern class can be represented by a simplicial 2-cocycle of $B$
having values 0 or 1. The simplicial set of circular permutations is
canonically a quotient of the simplicial set of all permutations $\pmb
S_*$ by a simplicial equivalence relation induced by \textit{right}
actions of cyclic subgroups. The simplicial set of all permutations
$\pmb S_*$ has the structure of a \textit{symmetric crossed simplicial
	group}. We have the simplicial map:
\begin{equation} \label{cyrc} \pmb S_* \xar{\ca} \pmb {SC}_* \end{equation}

We aim to prove the following:
\begin{theorem} \label{main} \begin{equation*} |\pmb{SC}_*| \approx
		K(\Z,2). \end{equation*} \end{theorem}

To the author's limited knowledge, $\pmb{SC}_*$ is the first simplicial
model of $K(\Z,2)$ with a finite number of simplices in every
dimension. This fact likely makes the simplicial set $\pmb{SC}_*$
interesting. The situation is somewhat related to the well-known topic
of triangulating $\C P^n$. See \cite{MorYo1991, ArMar1991} and the new
results in \cite{datta2024simplicial}. There are also interesting
computer experiments in \cite{Sergeraert2010}. The connections between
these results and our construction need further investigation. The
connection is probably through the minimal triangulation of the
tautological Hopf bundle $U(1)\xrightarrow{}S^{2n+1}\xrightarrow{}\C
P^n$.

Crossed simplicial group theory originated from pioneering works on
cyclic homology \cite{Tsygan1983} and \cite{Connes1983}. The idea
behind the proof of Theorem \ref{main} is to reference the remarkable
theorems on geometric realizations of crossed simplicial groups and
sets (\cite[Theorem 2.3]{Kras1987} \cite[Theorem 5.3, Lemma
5.6]{FL1991}, \cite[Theorem 7.1.4, Exercise 7.1.4 ]{Loday1998}). The
first mention of geometric realization for cyclic sets as
$U(1)$-spaces, and the main ingredient of the construction—the
\textit{geometric cyclic cosimplex}, or \textit{twisted
shuffle	product} $S_\cdot^1 \times_t \Delta^k \approx S_\cdot^1 \times
\Delta^k$ (here $S_\cdot ^1$ is a circle composed of one 1-simplex and
one point) is found in \cite[pp. 208-209]{Goodwillie1985} and further
extended in \cite[\S 2, Proposition 2.4]{DHK1985}, \cite[Theorem
3.4]{Jones1987}.

Theorem 1 immediately follows from an inspection of the constructions
in the above theorems. The arguments are geometrical.  As a result, we will see that the minimally
triangulated circle bundles over simplices described in \cite{MnevMin}
are nothing more than canonically \textit{order reoriented} twisted shuffle product $S_\cdot^1 \times_t \Delta^k$, and the map
(\ref{cyrc}) is the universal minimally triangulated circle bundle.

Section \ref{prelim}:
In this section, we discuss the basics of crossed simplicial group
theory for the case of $\pmb C_* \leq \pmb S_*$, recalling the left
crossed action of $\pmb C_*$ on $\pmb S_*$, left crossed cyclic orbits
in $\pmb S_*$, and their geometric realizations. Classical left crossed
cyclic orbits in $\pmb S_*$ \textit{do not} form a simplicial
equivalence relation and have no direct simplicial quotient.

Section \ref{proof}:
In this section, we will explain how to deal with the \textit{right}
action of $\pmb C_*$ on $\pmb S_*$ as opposed to the canonical situation
of the left action. Right orbits \textit{do} form a simplicial
equivalence relation. We obtain $\pmb{SC}_*$ as the simplicial quotient,
which is the set of right cyclic crossed orbits. After geometric
realization, $|\pmb{SC}_*|$ is the cellular structure on the set of right
orbits $|\pmb S_*|/|\pmb C_*|$. Here, $|\pmb S_*|\approx * $ is a
contractible Hausdorff topological group, and $|\pmb C_*| = U(1)$ is a
Lie group. Therefore, the quotient map $U(1)\xar{}|\pmb S_*|
\xar{|\ca|} |\pmb {SC_*}|$ is a $U(1)$-fibration, and
$|\pmb{SC_*}|\approx K(\Z,2)$, which concludes the proof of Theorem
\ref{main}.

\bigskip
Author is deeply grateful to Boris Tsygan and Andr\'e Henriques for pointing the author to the subject of crossed simplicial groups.

\section{Preliminaries} \label{prelim} 
The pair of the symmetric crossed simplicial group $\pmb S_*$ and its cyclic subgroup $\pmb C_*$, $\pmb C_* \leq \pmb S_*$ is specially discussed in \cite[6.1]{Loday1998}. 
 
\p{Simplicial notations.}
We denote $\pmb \Delta$ the category of finite linear orders $[n]=\{0,1,2,...,n\}$ and non-decreasing maps between them called operators. The category $\N$ is generated by ``cofaces" $\delta_i$  and ``codegeneracies" $\sigma_i$ :
$$ [n-1]\xar{\delta_i} [n] \xleftarrow{\sigma_j} [n+1], i,j = 0\ldots n $$
Cofaces $\delta_i$ are the only injective order preserving maps ``missing $i$" in the target.
Codegeneracies $\sigma_j$  are the only non-decreasing surjections ``hitting $j$ in the  target twice", i.e $\sigma_j(j)=\sigma_j(j+1)=j$.
Opposite category $\pmb \Delta^\op$ is generated by faces $d_i = \delta_i^\op$,  and degeneracies $s_i=\sigma_i^\op$.    
Simplicial set X  is a functor $\pmb \Delta^\op \xar{X} \pmb{Sets}$. Face and degeneracies goes to face and degeneracy maps which are again denoted $X_n \xar{d_i} X_{n-1}$ and $X_{n}\xar{s_i}X_{n+1}$.  The category of functors 
$\pmb \Delta^\op \xar{} \pmb{Sets}$ or ``presheaves" on $\pmb \Delta$ and  natural transformations of those (maps of simplicial sets) is denoted by $\widehat{\pmb \Delta}$.  
\p{Category $\pmb{\Delta G}$.}
Crossed simplicial groups and sets  are related to extension of $\pmb \Delta$ and   $\widehat{\pmb \Delta}$ by a correct adjoining of automorphism groups $\pmb G^\op_n $ to  $[n]$ in a such way that $\pmb G_n$ will act correctly  by automorphisms of sets $X_n$.    
\begin{definition} \label{Crossed} 
	\textbf{\cite[Definition  1.1]{FL1991}} \\
	A sequence of groups $\pmb G=\{\pmb G_n\}, n \geq 0$ is a crossed simplicial group if it is equipped with the following structure. There is a small category $\pmb {\Delta G}$, which is part of the structure, such that 
	\begin{description} 
		\item[(a)] the objects of $\pmb {\Delta G}$ are $[n]$, $n \geq 0$, 
		\item[(b)] $\pmb {\Delta G}$ contains $\pmb \Delta$ as a subcategory, 
		\item[(c)] $\on{Aut}_{\pmb {\Delta G}}([n])= \pmb G_n^\op$ (opposite group of $\pmb G_n$) , 
		\item[(d)] any morphism $[m]\xar{}[n]$ in $\pmb{\Delta G}$ can be uniquely written as a composite $\xi \cdot g$  where $ \xi \in \on{Hom}_{\pmb \Delta}([m],[n])$ and $f\in \pmb  G^\op_m$ (whence the notation $\pmb {\Delta G}$ ).
	\end{description}
	$$ 
	\begin{tikzcd}
		{[m]} \arrow[rd] \arrow[d, "g"', dotted] &       \\
		{[m]} \arrow[r, "\xi", dotted]          & {[n]}
	\end{tikzcd}
	$$ 
\end{definition}

\p{$\pmb{\Delta C} \subset \pmb{\Delta S}$, $\pmb C_* \leq \pmb S_* $.} \label{deltaS}
Here we follow \cite[6.1]{Loday1998}. We denote $\pmb S_n$ the group of permutations of $n+1$ (sic!) ordered elements $[n] =\{0,1,\ldots,n\}$, i.e. $g \in \pmb S_n$ is a one-to one map $[n]\xar{f}[n]$ represented as permutation $(f(0),\ldots,f(n))$. We have a commutative subgroup $\pmb C_n \leq \pmb S_n$ of cyclic permutations generated by the cycle $\tau = (n,0,1,\ldots,n-1)$.  We denote $\pmb S_*$, $\pmb C_*$ corresponding graded groups equipped with graded multiplication.  They are  equipped with simplicial structure interacting with multiplication   in a canonical ``crossed" way. For this we should pass to category $\pmb {\Delta S}$.      

\bigskip
The category $\pmb {\Delta S}$ is the category $\pmb \Delta$ enlarged by groups of arbitrary non-monotone automorphisms of ordered sets $[n]$ written as opposite symmetric group $\pmb S_n^\op$ (or $\pmb S_n$ acting from the right on $[n]$). Checking  and unwinding 
conditions of Defition \ref{Crossed} is subject of \cite[Theorem 6.1.4]{Loday1998}, see also \cite[Appendix A10 p. 191]{FT1987}. 

\begin{figure}[h!]
	\captionsetup{font=small,width=0.8\textwidth}	
	\begin{center}
		\includegraphics[width=4.0in]{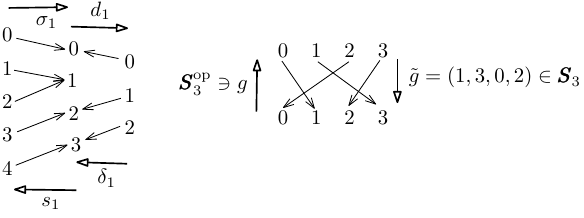} 
		\caption{ \label{permw} Wire diagrams of operators, permutations and their opposites.}
	\end{center}
\end{figure}

It is instructive to imagine both permutations and operators of $\pmb \Delta$ as ``wire diagrams" of maps between finite linear  orders  (Fig. \ref{permw}). 

\begin{figure}[h!]
	\captionsetup{font=small,width=0.8\textwidth}	
	\begin{center}
		\includegraphics[width=5.0in]{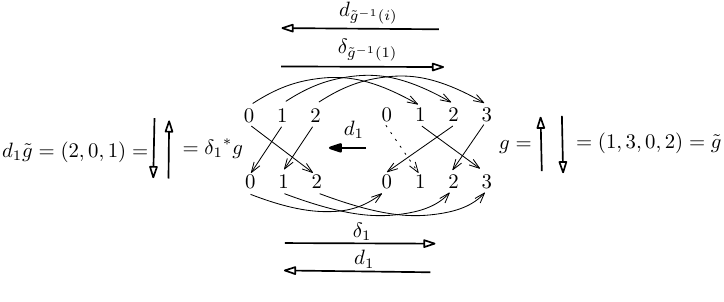}
`		\includegraphics[width=5.0in]{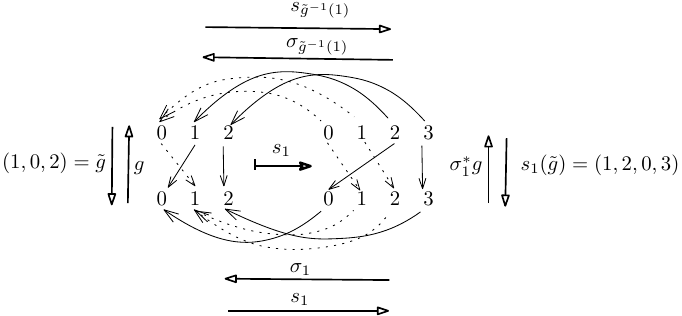} 
		\caption{ \label{b_d} In wire diagrams deletion $d_i$ corresponds to deletion of the arrow with target $i$ and degeneracy $s_i$ corresponds to parallel doubling of the arrow with target $i$. }
	\end{center}
\end{figure}
For an element $g \in \pmb S_n^\op$ there is associated  set map $[n]\xar{g}[n]$ with the same name $g(i) = \tilde g^{-1}(i)$, where $\tilde g \in \pmb S_n$ is the corresponding permutation.  

In the language of wire diagrams permutations and their duals in the opposite groups,  (co)boundaries  (co)degeneracies  communicate as depicted  in Fig. \ref{b_d}. Inspecting wire diagrams for (co)boundaries and (co)degeneracies we get that 
simplicial relations produces \textit{for any pair
 $( g,\xi)$, 
$ g \in \pmb S_n^\op, \xi \in \on{Hom}_{\pmb\Delta} ([m],[n])$
unique maps $\xi^* g,  g_*\xi$ such that \\
 (i) the following diagram is commutative:}   
\begin{equation} \label{decomp}
	\begin{tikzcd} 
		{[m]} \arrow[r, "{\xi \in \pmb{\Delta}}"] \arrow[d, "\pmb{S}^\op_m \ni \xi^*g "', dashed] & {[n]} \arrow[d, "{ g \in \pmb{S}_n^\op}"] \\
		{[m]} \arrow[r, " g_*\xi \in \pmb{\Delta}"', dotted]                  & {[n]}               
	\end{tikzcd}
\end{equation}
\textit{
and\\
(ii) restriction of $\xi^*g$ to each subset $\xi^{-1}(i), i=0...n$ preserves the order.}\\ 
The above  statement  is the subject of \cite[Lemma 6.1.5]{Loday1998}.

Thus we have a category $\pmb{ \Delta S}$ with objects - finite orders $[n]$ and morphisms - pairs $[m]\xar{(\xi,g)}[n]$, $\xi \in \on{Hom}_{\pmb \Delta}([m],[n])$, $g \in \pmb S^\op_m$.
Having another morphism $[k]\xar{(\phi,h)}[m]$
 the composition is defined by the rule 
 \begin{equation} \label{comp}
 	(\xi,g)\circ(\phi,h) = (\xi \circ g_*\phi,\phi^*\circ h )
 \end{equation} where the compositions of components  are in $\pmb \Delta$ and $\pmb S_k^\op$. The category $\pmb {\Delta S}$ satisfies requirements of Definition \ref{Crossed}. 
The opposite category $\pmb{\Delta S}^\op$
has decomposition of arrows opposite to (\ref{decomp}): 
\begin{equation} \label{codecomp}
	\begin{tikzcd}
		{[m]}                                             & {[n]} \arrow[l, "\alpha \in \pmb{\Delta}^{op}"']                                         \\
		{[m]} \arrow[u, "\pmb{S}_m \ni\alpha_*f", dashed] & {[n]} \arrow[u, "f \in \pmb{S}_n"'] \arrow[l, "f^*\alpha \in \pmb{\Delta}^{op}", dashed]
	\end{tikzcd}
\end{equation}
When in (\ref{codecomp}) we set $\alpha = d_i = \delta_i^\op$, $m=n-1$, we have 
$f^*d_i = d_{f^{-1}(i)}$ and we denote $(d_i)_* f$  by $d_i f$. When in (\ref{codecomp}) we choose $\alpha = s_i = \sigma_i^\op$, $m=n+1$, we have  
$f^*s_i = s_{f^{-1}(i)}$ and we denote $(s_i)_* f$ by $s_i f$. 
Applying opposite to composition rule (\ref{comp}) we got that 
for $f,h \in \pmb S_n$
\begin{equation}
	\begin{aligned} \label{dmult}
 d_i (h\circ f) & = d_i h\circ d_{h^{-1}(i)}f \\
  s_i (h\circ f) & = s_i h\circ d_{h^{-1}(i)}f		
	\end{aligned}
\end{equation}   
Replacing symmetric groups $\pmb S_n$ by cyclic subgroups $\pmb C_n$ we obtain subcategory $\pmb {\Delta C} \subset \pmb{\Delta S}$.

Now crossed simplicial group $\pmb S_*$ can be canonically identified with representable $\pmb{Sets}$-valued Yoneda presheaf   \begin{equation} \label{ys}
	 \pmb{Y}_{\pmb{\Delta S}}([0])=\on{Hom}_{\pmb{\Delta S}}(-,[0]) = \pmb S_*
\end{equation} By decomposition rules it is simplicial set structure on graded set of groups  $\pmb S_*$ with boundaries and degeneracies defined by (\ref{codecomp}) and communicating with multiplication in a ``crossed" way by rules (\ref{dmult}) (see \cite[Proposition 1.7]{FL1991}). On $\pmb C_*$ we have induced structure
\begin{equation} \label{yc}
	 \pmb{Y}_{\pmb{\Delta C}}([0])=\on{Hom}_{\pmb{\Delta C}}(-,[0]) = \pmb C_* \leq \pmb S_*
\end{equation}
 thus the pair $\pmb C_* \leq \pmb S_*$ is defined.
\p{$\pmb C_* \leq \pmb S_*$ in terms of permutations.}
Here we  rephrase the resulting from canonical $\pmb{\Delta S}$-construction \S(\ref{deltaS})  structure of $\pmb S_*$ in terms of permutations. So,  denote $\pmb S_n$ the group of permutations of $n+1$ ordered elements $[n] =\{0,1,\ldots,n\}$, i.e. $f \in \pmb S_n$ is a one-to one map $[n]\xar{f}[n]$ represented as permutation $(f(0),\ldots, f(n))$. The graded set of permutations $\pmb S_* = \pmb S_0,\pmb S_1 \ldots$ forms a simplicial set. The $i$-th boundary map $\pmb S_n \xar{d_i}\pmb S_{n-1}, i=0,...,n $ is deleting $i$-th element of permutation and reordering other elements monotonically, i.e. elements from $0$ to $i-1$ preserves their numbers. Elements from $i+1$ to $n$ got the numbers $i\ldots n-1$ (see (\ref{permbd})). The $i$-th degeneracy $\mf S_n \xar{s_i} \mf S_{n+1}, i=0,...,n$ inserts element with number $i+1$ \textit{next} to the element $i$ and  reorders other elements monotonically. Elements from $0$ to $i$ preserves numbers and the old elements $i+1 ... n$ of the permutation  got shifted by one numbers $i+2...n+1$ correspondently.
\begin{equation} \label{permbd}
	\begin{aligned}
		(d_i f)(j) = & 
		\begin{cases}
			f(j) & \text{ if } j=0\ldots i-1 \\
			f(j-1) & \text { if } j = i \ldots n \\
		\end{cases} \\	
		(s_if)(j) = & \begin{cases}
			f(j) & \text{ if } f(j)=0,...,i\\
			i+1    & \text{ if } j = f^{-1}(i)+1 \\
			f(j)+1 & \text { if } f(j)=i+1,...,n 
		\end{cases} \\
	\end{aligned}             
\end{equation}    
Additionally in $\pmb S_*$  we have crossed multiplication  $(f_n, g_n) \mapsto f_ng_n$ communicating with boundaries and degeneracies by rules (\ref{dmult}). This crossed multiplication  will became canonically  functorial in \S(\ref{crossedprod})   
We have a crossed simplicial  subgroup $\pmb C_n \leq \pmb S_n$ of cyclic permutations generated by cycles $\tau_n = (n,0,1,\ldots,n-1)$.

\p{ $\pmb{S}_* \xar{\ca} \pmb{SC}_*$} \label{SC}
We recall the simplicial  map  $\pmb{S}_* \xar{\ca} \pmb{SC}_*$ from \cite{MnevMin}. The group $\pmb C_n$ acts from the \textit{right} on 
permutations by shifts $f_n \tau_n = (f_n(n),f_n(0),\ldots f_n(n-1))$. The orbits of the right action of $\pmb C_n$ on $\pmb S_n$ are numbered by \textit{circular permutations}, i.e. oriented circular necklaces with $n+1$ beads coloured by $[n]$. We denote this set of right orbits or $n+1$ circular permutations $\pmb S_n / \pmb C_n$ by $\pmb{SC}_n $. 
The rules (\ref{permbd}) induces  the simplicial set structure on the graded set of circular permutations $\pmb {SC}_* = \pmb {SC}_0, \pmb {SC}_1 \ldots$ -- we can delete a bead $i$ (this provides $d_i$ ) and we can insert a bead $i+1$ right after the bead $i$ since luckily the relation   ``right after" exist in circular order (this provides $s_i$).  Thus we got simplicial set of circular permutations $\pmb{SC}_*$ together with the  simplicial factor-map $\pmb{S}_* \xar{\ca} \pmb{SC}_*$ sending a permutation to its right  cyclic orbit.

\p{Simplicial, cyclic and symmetric sets, base change adjacency and left crossed cyclic orbit of a permutation.}\mbox{} \label{basech}

Symmetric or cyclic set is a $\pmb{Sets}$-valued presheaf on $\pmb{\Delta S}$ or $\pmb{\Delta C}$.
In the  following we use  $\pmb G$ for definitions and  statements which are equivalent for $\pmb S$  and $\pmb C$.
For example $\pmb G_*$ (\ref{ys}),(\ref{yc}) is a $\pmb G$-set.

\bigskip
The important point for us is that due to embedding $\pmb{\Delta C} \subset \pmb{\Delta S}$ we got that canonically $\pmb S_*$ is a $\pmb C$-set.

\bigskip
The categories of of $\pmb G$-sets with morphisms - natural transformations are denoted by $\widehat{\pmb{\Delta G}}$. By construction these are simplicial sets $X$ with fixed left actions of groups $\pmb G_n$ by automorphisms of $X_n$.   The action are explicitly described  by ``base change adjunction". 
\subsubsection{Adjunction data.}
We recall (see \cite[Chapter X]{MacLane98}) that adjunction $\la F,G,\varphi\ra$ 
between two small categories $ A, B$ is a pair of functors 
$A\xar{F}B, B\xar{G}A$ and  bifunctorial isomorphism of $\on{Hom}$ sets $B(F(X),Y)\xar{\varphi} A(X,G(Y))$, where $X$ is running over $A$ and $Y$ over $B$. Functor $F$  called left adjoint to $G$, $G$ called right adjoint to $F$, and adjunction sometimes denoted by $F \dashv G$. Adjunction defines and is defined by ``monad of adjunction": a natural transformation of $A$-endofunctors $\on{Id}_A \xar{\iota} GF$ called ``unit of adjunction" and a natural transformation $B$-endofunctors $FG\xar{\varepsilon} \on {Id}_B $ called counit of adjunction
satisfying ``triangular identities"    
\begin{equation} \label{trieq}
	\begin{tikzcd}
		F \arrow[r, "F \cdot \iota"] \arrow[rd, "id"] & F\circ G \circ F \arrow[d, "\varepsilon \cdot F"] \\
		& F                                                                                 
	\end{tikzcd}\begin{tikzcd}
		G \arrow[r, "\iota \cdot G "] \arrow[rd, "id"] & G\circ F \circ G \arrow[d, " G \cdot \varepsilon"] \\
		& G                                                                                  
	\end{tikzcd}
\end{equation} 

Embedding of (skeletal) categories $\pmb \Delta \subset \pmb {\Delta G}$ creates embedding of those duals $\pmb \Delta^\op \xar{\mc P} \pmb {\Delta G}^\op$. On presheaves we got forgetful functor 
$\widehat{\pmb \Delta} \xleftarrow{ \ol{*}} \widehat{\pmb {\Delta G}}$ making simplicial set $\ol Y$ from $\pmb G$-set $Y$. Functor $\ol {*}$ has left adjoint which we denote\footnote{The functor is denoted by $G$ in \cite{Kras1987} and $F$ in \cite{FL1991}} $\pmb G_* \times_t *$:  $(\pmb G_* \times_t *) \dashv \ol{*}$. The left adjoint is computed as pointwise left Kan extension of simplicial set $X$ along $\mc P$ \cite[X.3 Theorem 1]{MacLane98}. This is a specially simple situation of \textit{``base change adjunction"}.
\subsubsection{Crossed left action.}
The left Kan extension of $X$ along $\mc P$ produces the following element-wise  formulas for simplicial and 
$\pmb G_*$-structure on  $\pmb G_*\times_t X$: 
\begin{equation}
	\begin{aligned}
		(\pmb{G}_*\times_t  X)_n & =\{(h_n, x_n)\}_{h_n \in \pmb G_n, x_n \in X_n}\\
		d_i (h_n, x_n) & = (d_i h_n, d_{h^{-1}(i)}x_n) \\	
		s_i (h_n, x_n) & = (s_i h_n, s_{h^{-1}(i)}x_n)\\
		f_n \cdot (h_n, x_n) & = (f_nh_n, x_n)	
	\end{aligned}
\end{equation}
The adjunction $(\pmb G_* \times_t *) \dashv \ol *$ defines monad  
with the unit
$$\on {id}_{\widehat{\pmb \Delta}} \xar{\iota} \ol{(\pmb G_* \times_t * )}$$ and counit $$(\pmb G_* \times_t \ol *) \xar{ev} \on {id}_{\widehat{\pmb \Delta G}}$$ satisfying ``triangular identities" (\ref{trieq}). 
The unit of the adjunction is computed on elements as follows.
For an element $x_n \in X_n$ of simplicial set $X$ we got
$$\iota (x_n) = (1_{\pmb G_n}, x_n ) $$ Counit of the adjunction defines  \textit{the crossed  left action} of $\pmb G_*$ on $\pmb G$-set $Y$, namely for $y_n \in Y_n$ we got   
$$ev (g_n, y_n) = g_n\cdot y_n$$
 Triangular identities  (\ref{trieq}) of the monad ensures that the action is correct action of $\pmb G_*$ in a crossed way:
$$\begin{tikzcd}
	\pmb G_* \times_t X \arrow[rr, "{(h,x)\mapsto (h, (1,x))}"] \arrow[rrd, "\on{id}"] &  & \pmb G_* \times_t \ol{(\pmb G_* \times_t X)} \arrow[d, "{(h,(f,x)) \mapsto (hf, x))}"] \\
	&  & \pmb G_* \times_t X                                                          
\end{tikzcd}
\begin{tikzcd}
	X \arrow[r, "{x \mapsto (1,x)}"] \arrow[rd, "\on{id}"] & \ol {\pmb G_* \times_t X} \arrow[d, "{(h,x) \mapsto h\cdot x}"] \\
	& X                                                      
\end{tikzcd}$$

\subsubsection{Crossed product.} \label{crossedprod}
If $X=\pmb G_*$ then the  counit  
$$\pmb G_*\times_t \pmb G_* \xar{ev} \pmb G_*$$ represents ``crossed product" in crossed simplicial group $\pmb G_*$. 
\subsubsection{Yoneda Lemma and left crossed orbits.} \label{leftorb}
Categories of presheaves $\widehat {\pmb \Delta} $, $\widehat {\pmb \Delta G} $ has representable (Yoneda) objects -- cosimplices $$\pmb{\Delta}[n] = \pmb{\Delta}(-,[n])=\pmb{\Delta}^\op([n],-)$$ and $\pmb{G}$-cosimplices 
$$\pmb{\Delta G}[n] = \pmb{\Delta G}(-,[n])=\pmb{\Delta G}^\op([n],-)$$ 
There is the key isomorphism $\pmb {\Delta G}[n] \approx \pmb G_* \times_t \pmb{\Delta}[n]$ (\cite[Exersise 4.5]{FL1991}). 
(Co)Yoneda Lemma states that every presheaf is canonically colimit of representables. For simplicial set $X$ and $x_n \in X_n$ this creates colimit cone structure map $$\pmb{\Delta}[n]\xar{\pmb y_{\pmb \Delta}(x_n)} X$$ sending $\on{id}[n]$ to $x_n$.  Analogously for $\pmb G$-set $Y$ and $y_n \in Y_n$ this creates colimit cone structure map 
$$\pmb{\Delta G}[n]\approx \pmb G_*\times_t \pmb \Delta[n]\xar{\pmb y_{\pmb \Delta G}(y_n)} Y  $$ sending $(1_{\pmb G_n},\on{id}[n])$ to $y_n$.
Relation between unit-counit of adjunction and bifunctorial isomorphism  of $\on{Hom}$-sets
$$\widehat{\pmb{\Delta G}}(\pmb G_* \times_t X, Y) \xar{\varphi} \widehat{\pmb \Delta} (X, \ol Y)$$
connects the two Yoneda maps. For a $\pmb G$-set Y and element $y_n \in Y_n$ the isomorphism $\varphi$ sends $\pmb G_* \times_t \pmb \Delta[n] \xar{\pmb y_{\pmb {\Delta G}}(y_n) }Y$ to $\pmb \Delta[n] \xar{\pmb y_{\pmb \Delta} (\ol y_n)} \ol Y$. In the inverse direction $\varphi$ sends $\pmb y(\ol y_n)$ to $\pmb y_{\pmb \Delta G}(y_n)$ by the following commutative diagram 
\begin{equation}
\begin{tikzcd} \label{twoyo}
	\pmb G_*\times_t \overline Y \arrow[r, "ev"]                                                                                             & Y \\
	{\pmb G_*\times_t\pmb \Delta[n]} \arrow[u, "\pmb G_*\times_t \pmb y_{\pmb \Delta}(\ol y_n)"] \arrow[ru, "\pmb y_{\pmb \Delta G}(y_n)"'] &  
\end{tikzcd}
\end{equation}
 The Yoneda map $\pmb y_{\pmb{\Delta G}}(y_n)$ and its image in $\pmb G$-set Y we call \textit{left crossed $\pmb G$-orbit of $y_n \in Y_n$}.   
\p{Geometric realization.} \label{geomr}
Here we in situation of \cite[Theorem 2.3]{Kras1987}, \cite[Teorem 5.3]{FL1991}.
Geometric realization $|X| $ of a $\pmb G$-set $X$ is the geometric realization $|\overline X|$ of the underground simplicial set.     
The core of geometric  realization theorems states that there  is a canonical functorial homeomorphism\footnote{ The homeomorphism $\Psi$  is the 
	homeomorphism   ${\Phi(X)}^{-1}$ in \cite[Theorem 2.3]{Kras1987} 
	and $(p_1,p_2)^{-1}$ in \cite[Teorem 5.3]{FL1991}.}
$$|\pmb G_*|\times |X| \xar{\Psi} |\pmb G_* \times_t X|$$
such that in induced from geometric realization  metric the composite 
$$|\pmb G_*|\times |\pmb G_*| \xar{\Psi} |\pmb G_* \times_t \pmb G_*| \xar{|ev|} |\pmb G_*|$$ is a topological group.
For any $\pmb G$-set $Y$ the composite 
$$|\pmb G_*|\times |Y| \xar{\Psi} |\pmb G_* \times Y| \xar{|ev|} |Y|$$ 
makes $|X|$ left topological $|\pmb G_*|$-space.

\bigskip
In our situations $|\pmb C_*|$ is an oriented circle $S^1_\cdot$ made from one vertex and one non-degenerate 1-simplex (and oriented by its orientation). In induced metric the composite map
 $$|\pmb C_*|\times |\pmb C_*| \xar{\Psi} |\pmb C_* \times_t \pmb C_*| \xar{|ev|} |\pmb C_*|$$ is exactly $\R/\Z \approx U(1)$ group structure on $S_\cdot^1$ with the unit in the vertex of $S_\cdot^1$.
 $$|\pmb S_*|\times |\pmb S_*| \xar{\Psi} |\pmb S_* \times_t \pmb S_*| \xar{|ev|} |\pmb S_*|$$ is a contractible  topological group (\cite[Example 6]{FL1991})
 and since $\pmb S_*$ is a $\pmb C$-set the induced composed 
 map
 $$U(1)\times |\pmb S_*| \xar{\Psi} |\pmb C_* \times_t \pmb S_*| \xar{|ev|} |\pmb S_*|$$ is a free left action of Lie subgroup $U(1)\leq |\pmb S_*|$  on Hausdorff contractible space  $|\pmb S_*|$.
 
 We are specially interested in orbits of the action. 
 Cyclic cosimplex $\pmb y_{\pmb{\Delta C}}=\pmb \Delta \pmb C_*[n] = \pmb C_* \times_t \pmb \Delta[n]$ is a cyclic set. 
 Its geometric realization has cellular structure of ``twisted shuffle product" $S_\cdot^1\times_t \Delta^n$ (\cite[pp. 208-209]{Goodwillie1985},\cite[\S 2, Proposition 2.4]{DHK1985}, \cite[Theorem 3.4]{Jones1987}). Applying geometric realization to (\ref{twoyo}) we get a comutative diagram of spaces:
  \begin{equation}
 \begin{tikzcd}
 	U(1)\times{|\pmb S_*|}  \arrow[r, "\Psi"]                                                              & {|\pmb C_*\times_t \pmb S_*|} \arrow[r, "{|ev|}"]                                                                                                                     & {|\pmb S_*|} \\
 	U(1)\times \Delta^n \arrow[r, "\Psi"] \arrow[u, "id \times |\pmb y_{\pmb \Delta}(g_n)|" description] & {|\pmb C_*\times_t\pmb \Delta[n]|} \arrow[u, "|\pmb C_*\times_t \pmb y_{\pmb \Delta}(g_n)|\;\;\;\;\;" description] \arrow[ru, "{|\pmb y_{\pmb {\Delta C}}(g_n)|}"'] &           
 \end{tikzcd}
 \end{equation}
 Thus the composit map
 \begin{equation}
U(1)\times \Delta^n \xar{\Psi}S^1\times_t\Delta^n=|\pmb C_* \times_t \pmb{\Delta}[n]| \xar{|\pmb y_{\pmb{\Delta C}}(g_n)|} |\pmb S_*|
 \end{equation} 
 provides as its image a continuous trivial   family of left $U(1)$ orbits on $|\pmb S_*|$ parametrized by cell which is the image of characteristic map $\Delta^n\xar{|\pmb y_{\pmb \Delta}|}|\pmb S_*|$. In this way in geometric realization the crossed left $\pmb C_*$-orbit became a trivial family of true $U(1)$ left orbits.

\section{Proof of Theorem \ref{main}.} \label{proof}

\p{Left and right orbit spaces of topological subgroup.} \label{lr}
\label{left-right}

For a group \(G\) and a subgroup \(H \leq G\), there are left and
right actions of \(H\) on \(G\), and these actions are free. The left
action of \(H\) on \(G\) creates the set of left orbits \(G \backslash
H = \{Hg\}_{g \in G}\), while the right action of \(H\) on \(G\)
creates the set of right orbits \(G / H = \{gH\}_{g \in G}\). The
group \(G\) acts from the right on \(G \backslash H\) and from the
left on \(G / H\), with stabilizer \(H\).

The involution
\begin{equation} \label{Upsilon}
	G \xrightarrow{\upsilon} G : \upsilon(g)=g^{-1}
\end{equation}
switches between left and right \(H\)-orbits of \(g\) and \(g^{-1}\).
It is standard to define the \textit{opposite group} \(G^\op\) with
the same elements as \(G\) but with multiplication \(g_1 * g_2 = g_2
g_1\). Then \(\upsilon\) is a group isomorphism \(G
\xrightarrow{\upsilon} G^\op\) that maps left \(H\) orbits in \(G\) to
left \(H^\op\) orbits in \(G^\op\) (which were right \(H\)-orbits in
\(G\)), thereby inducing a one-to-one correspondence
\begin{equation} \label{upsilon}
	G \backslash H \overset{\tilde \upsilon}{\approx} G / H
\end{equation}
between the sets of right and left \(H\)-orbits in \(G\).

In the topological category, where \(H\) and \(G\) are topological
groups, the sets \(G \backslash H\) and \(G / H\) become orbit spaces
with quotient topology, and \(\tilde \upsilon\) in (\ref{upsilon}) is
a homeomorphism between left and right orbit spaces. In good
situations, for example, if \(H\) is a Lie group and \(G\) is
Hausdorff, the map
\begin{equation} \label{top_factor}
	H \xrightarrow{} G \xrightarrow{} G \backslash H \;\; [G / H]
\end{equation}
is locally trivial with fiber \(H\) (see \cite[4.1 on page
315]{Palais1961}). Therefore, \(G \backslash H \;\; [G / H]\) is
Hausdorff (see \cite[Theorem 31.2 (a) on page 196]{Munkres2000}), and
(\ref{top_factor}) is a principal Serre fibration (see \cite[Satz
5.14]{DieckPuppe}).

In the simplicial category, where \(H\) and \(G\) are simplicial
groups, \(G \backslash H\) and \(G / H\) have the structure of
simplicial sets, and \(H \xrightarrow{} G \xrightarrow{} G \backslash
H \;\; [G / H]\) is a principal Kan fibration (see \cite[Definition
18.1, Lemma 18.2]{May1968}).

\p{Left vs. right crossed action problem.}

In the crossed-simplicial setting, a cyclic crossed simplicial group
\(\pmb C_*\leq \pmb S_*\) is a subgroup of a symmetric
crossed-simplicial group, acting on \(\pmb S_*\) from the
\textit{left} in a twisted manner. This twist disappears in geometric
realization (see \S(\ref{geomr})). The geometric realization theorems
for crossed simplicial groups imply that \(U(1) \approx |\pmb{C}_*|
\leq |\pmb{S}_*|\), where \(|\pmb S_*|\) is a contractible topological
group. This provides, according to \S(\ref{left-right}), a principal
\(U(1)\) fibration
\begin{equation} \label{left}
	U(1) \approx |\pmb C_*| \xrightarrow{} |\pmb S| \xrightarrow{}
	|\pmb S_*| \backslash |\pmb C_*| \approx K(\Z,2)
\end{equation}
However, for the \textit{left} crossed action of \(\pmb C_*\) on
\(\pmb S_*\), there is \textit{no simplicial structure} on the left
orbit set, since the left crossed orbits (\S(\ref{leftorb}))
\textit{are not} classes of simplicial equivalence relations.
Therefore, something like \(\pmb S_* \backslash \pmb C_*\) does not
exist simplicially.

On the other hand, for a crossed simplicial group, the opposite group
is not well-defined as a crossed simplicial group, so the switch
between left and right actions is not entirely trivial in the crossed
simplicial setting. We need the \textit{right} action to handle \(\pmb
S_* \xrightarrow{\ca} \pmb{SC}_*\) (\S (\ref{SC})).

\p{Universal order reorientation $\Upsilon$ of simplicial sets.}

\begin{figure}[h!]
	\captionsetup{font=small, width=0.8\textwidth}
	\begin{center}
		\includegraphics[width=4.0in]{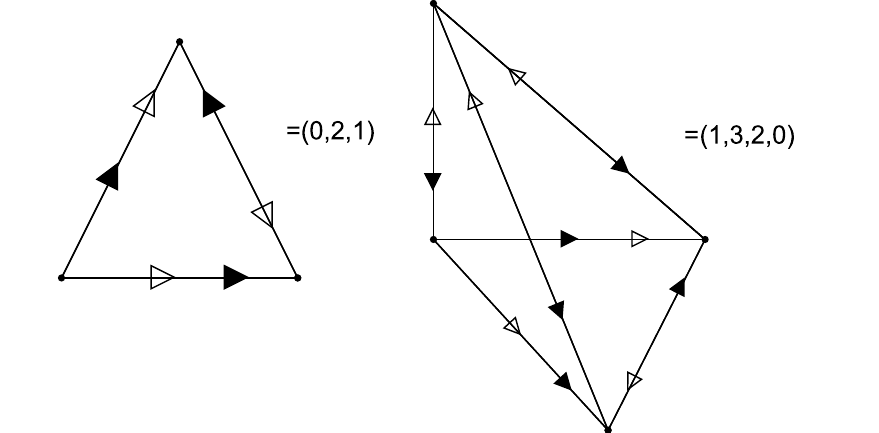}
		\caption{ \label{twoorders} Permutation as a double ordered
			simplex. Here
			\(\vartriangleright\) denotes the source order, and
			\(\blacktriangleright\) denotes the target order.}
	\end{center}
\end{figure}

A permutation can be identified with a simplex \(\Delta^n\) having
\textit{two} total orientation orders on vertices: source order and
target order (see Fig.\ref{twoorders}).

If one has a finite simplicial complex \(K\), a local order
orientation of its simplices can be obtained by fixing a total order
of its vertices, resulting in a semi-simplicial (or “\(\Delta\)”) set,
denoted \(K_{\text{source}}\). This is a standard method. Not all
local order orientations of \(K\) can be obtained this way. A
different total order on the vertices will yield different local
orders on simplices, resulting in a semi-simplicial set structure
\(K_{\text{target}}\) on the same complex \(K\). The two “source” and
“target” orders will produce source and target orders on the vertices
of every simplex of \(K\), both comparable with boundaries. Hence,
every \(n\)-simplex \(x_n\) of \(K\) has a permutation \(p(x_n)\),
providing a simplicial map \(K_{\text{target}} \xrightarrow{p}
\ul{\pmb S}_*\) (where \(\ul{\pmb S}_*\) is the semi-simplicial set
obtained from \(\pmb S_*\) by forgetting degeneracies). There is a
non-simplicial involution \(\ul{\pmb S}_* \xrightarrow{\Upsilon}
\ul{\pmb S}_*\) that switches between the source and target
permutations, i.e., sending a permutation \(f\) to \(f^{-1}\). The
involution induces a simplicial map \(K^{\Upsilon p} =
K_{\text{source}} \xrightarrow{p^{-1}} \ul{\pmb S}_*\). Thus, together
with the involution \(\Upsilon\), the semi-simplicial set \(\ul{\pmb
	S}_*\) “represents” a representable functor of the “double orientation
ordering” of a semi-simplicial set, along with the operation of
order reorientation.

The same involution acts on \(\pmb S_*\), respecting degeneracies. By
inspecting wire diagrams (see Fig \ref{b_d}), we can see that the
diagram of permutation \(f_n^{-1}\) is obtained from the diagram of
\(f_n\) by reversing the direction of arrows. In this process,
boundaries map to boundaries, and degeneracies to degeneracies in a
canonical but non-simplicial way. Thus, \(\Upsilon\) is a
non-simplicial automorphism of the simplicial set \(\pmb S_*\),
sending
\(f_n\) to
\[
\Upsilon(f_n) = f_n^{-1},
\]
boundary \(d_i f_n\) to boundary
\[
\Upsilon(d_i f_n) = d_{f_n^{-1}(i)} f_n^{-1},
\]
degeneracy \(s_i f_n\) to degeneracy
\[
\Upsilon(s_i f_n) = s_{f_n^{-1}(i)} f_n^{-1}
\]
providing a coordinate change on geometric realization. Also, we have
the left-right multiplication involution
\[
\Upsilon(f_n h_n) = h_n^{-1} f_n^{-1}.
\]

If one has a simplicial map \(X \xrightarrow{a} \pmb S_n\), this means
that simplices of \(X\) are decorated by permutations in a way
compatible with boundaries and degeneracies. We can
reorient simplices by changing the source and target orders, i.e., by
using the non-simplicial map \(X \xrightarrow{a} \pmb S_*
\xrightarrow{\Upsilon} \pmb S_*\), resulting in a new simplicial set
\(X^{\Upsilon a}\) on the same set of simplices, with canonically
homeomorphic geometric realization. Together with the involution \(\Upsilon\), the simplicial
set \(\pmb S_* \) represents  functor of double
orientation ordering of simplicial sets, along with the operation of
order reorientation.

\p{Right crossed \(\pmb C_*\)-orbits in \(\pmb S_*\).} \label{rightorb}

We don't know exactly what the opposite of a crossed simplicial group
is (since it is not a crossed simplicial group), but we can define a
\textit{right \(\pmb C_*\)-orbit} of \(g \in \pmb S_n\). For this, we
define a simplicial set denoted by \(E(\ca g)\). We follow notations
(\ref{codecomp}) for \(\pmb{\Delta S}^\op\). Define
\[
E(\ca g)_m = \{([n] \xrightarrow{\alpha} [m], \alpha_* g \cdot h) \mid
h \in \pmb C_m\}
\]
\[
d_i ([n] \xrightarrow{\alpha} [m], \alpha_* g \cdot h) = (d_i \alpha,
d_i(\alpha_* g \cdot h))
\]
\[
s_i ([n] \xrightarrow{\alpha} [m], \alpha_* g \cdot h) = (s_i \alpha,
s_i(\alpha_* g \cdot h))
\]
It has simplicial projections
\[
E(\ca g) \xrightarrow{q_1} \pmb \Delta [n]
\]
\[
E(\ca g) \xrightarrow{q_2 (g)} \pmb S_*
\]
The right-\(\pmb C_*\) orbit of \(g\) is by definition the image of
\(q_2(g)\) in \(\pmb S_n\).

Tautological computations provide the following lemma:
\begin{lemma}
	\mbox{}\\
	(i) Let \(\ca g_n \in \pmb{SC}_n\) and \(\Delta[n] \xrightarrow{\pmb
		y_{\pmb \Delta} (\ca g_n)} \pmb{SC}\) be the Yoneda simplex of \(\ca
	g_n\) in \(\pmb{SC}\). Then \(q_1\) and \(q_2(g_n)\) are the
	components of the pullback diagram
\begin{equation*} \label{Ecirc}
		\begin{tikzcd}
		E(\circlearrowright g_n) \arrow[d, "q_1"', dashed] \arrow[r,
		"q_2(g_n)", dashed] & \pmb S_* \arrow[d, "\circlearrowright"] \\
		\pmb \Delta[n] \arrow[r, "y(\circlearrowright g_n)"] & \pmb{SC}
	\end{tikzcd}
\end{equation*}
	(ii) \(\Upsilon (\pmb y_{\pmb{\Delta C}}(g_n)) = q_2(g_n^{-1})\), $E(\ca g^{-1}_n) = (\pmb C_* \times_t \pmb \Delta [n])^{\Upsilon \pmb y_{\pmb{\Delta C}}(g_n)}$.
This means that order orientation involution $\Upsilon$ turns 
left crossed cyclic orbit of permutation 
$\pmb C_*\times_t \pmb{\Delta}[n]\xar{\pmb y_{\pmb{\Delta C}}(g_n)} \pmb S_*$ into  pullback	
\begin{equation*} \label{Ecircc}
	\begin{tikzcd}
		E(\circlearrowright g^{-1}_n) \arrow[d, "q_1"', dashed] \arrow[r,
		"q_2(g^{-1}_n)", dashed] & \pmb S_* \arrow[d, "\circlearrowright"] \\
		\pmb \Delta[n] \arrow[r, "y(\circlearrowright g^{-1}_n)"] & \pmb{SC}
	\end{tikzcd}
\end{equation*}	
\end{lemma}

It follows that crossed right \(\pmb C_*\) orbits \textit{form} a
simplicial equivalence relations on \(\pmb S_*\) (unlike the left
orbits). Its factor set is \(\pmb{SC_*} \approx \pmb S_*/\pmb C_*\).
The space \(|E (\ca G)| \xrightarrow{p_1} \Delta^n\) is a minimally
triangulated circle bundle associated with \(\ca g_n\). It is just the
\(\Upsilon\)-reoriented geometric twisted shuffle product
\[
S_\cdot^1 \times_t \Delta^n = |\pmb C_* \times_t \pmb \Delta[n]|
\approx U(1) \times \Delta^n.
\]

\p{Order reorientation \(|\Upsilon|\) on geometric realization \(|\pmb
	S_*|\) is the canonical group involution \(\upsilon\).}

It follows from classical constructions (see \cite[the map \(\chi\) in
the proof of Theorem 2.3 on page 52]{Kras1987}) that the geometric
realization $|\Upsilon|$ of the order reorientation involution \(\Upsilon\) is
exactly the involution \(\upsilon\) ((\ref{Upsilon}) \S(\ref{lr})):
\begin{equation}
	\begin{tikzcd}
		\pmb C_* \leq \pmb S_* \arrow[r, "\Upsilon"] \arrow[d, "|*|"'] &
		\pmb C_* \leq \pmb S_* \arrow[d, "|*|"] \\
		{|\pmb C_*|} \leq {|\pmb S_*|} \arrow[r, "\upsilon"] & {|\pmb C_*|} \leq
		{|\pmb S_*|}
	\end{tikzcd}
\end{equation}
extending the chain (\ref{left}) by
\begin{equation} \label{right}
	|\pmb S_*| \backslash |\pmb C_*| \overset{\tilde \upsilon
	}{\approx} |\pmb S_*| / |\pmb C_*| \approx |\pmb {SC}_*| \approx
	K(\Z,2)
\end{equation}
This completes the proof of Theorem \ref{main}.

\def\cprime{$'$} \def\cprime{$'$}


\begin{thebibliography}{{May}68}
	
	\bibitem[AM91]{ArMar1991}
	Pierre Arnoux and Alexis Marin.
	\newblock The {K{\"u}hnel} triangulation of the complex projective plane from
	the view point of complex crystallography. {II}.
	\newblock {\em Mem. Fac. Sci., Kyushu Univ., Ser. A}, 45(2):167--244, 1991.
	
	\bibitem[Con83]{Connes1983}
	Alain Connes.
	\newblock Cohomologie cyclique et foncteurs {{\(Ext^ n\)}}.
	\newblock {\em C. R. Acad. Sci., Paris, S{\'e}r. I}, 296:953--958, 1983.
	
	\bibitem[DHK85]{DHK1985}
	W.G. {Dwyer}, M.J. {Hopkins}, and D.M. {Kan}.
	\newblock {The homotopy theory of cyclic sets.}
	\newblock {\em {Trans. Am. Math. Soc.}}, 291:281--289, 1985.
	
	\bibitem[DS24]{datta2024simplicial}
	Basudeb Datta and Jonathan Spreer.
	\newblock Simplicial cell decompositions of $\mathbb{CP}^{\hspace{.3mm}n}$,
	2024.
	
	\bibitem[FL91]{FL1991}
	Zbigniew {Fiedorowicz} and Jean-Louis {Loday}.
	\newblock {Crossed simplicial groups and their associated homology.}
	\newblock {\em {Trans. Am. Math. Soc.}}, 326(1):57--87, 1991.
	
	\bibitem[FT87]{FT1987}
	B.~L. {Feigin} and B.~L. {Tsygan}.
	\newblock {Additive \(K\)-theory.}
	\newblock {\(K\)-theory, arithmetic and geometry, Semin., Moscow Univ. 1984-86,
		Lect. Notes Math. 1289, 67-209 (1987).}, 1987.
	
	\bibitem[Goo85]{Goodwillie1985}
	Thomas~G. Goodwillie.
	\newblock Cyclic homology, derivations, and the free loopspace.
	\newblock {\em Topology}, 24:187--215, 1985.
	
	\bibitem[{Jon}87]{Jones1987}
	John D.~S. {Jones}.
	\newblock {Cyclic homology and equivariant homology.}
	\newblock {\em {Invent. Math.}}, 87:403--423, 1987.
	
	\bibitem[{Kra}87]{Kras1987}
	R.~{Krasauskas}.
	\newblock {Skew-simplicial groups.}
	\newblock {\em {Lith. Math. J.}}, 27(1):47--54, 1987.
	
	\bibitem[{Lod}98]{Loday1998}
	Jean-Louis {Loday}.
	\newblock {\em {Cyclic homology. 2nd ed.}}
	\newblock Berlin: Springer, 2nd ed. edition, 1998.
	
	\bibitem[{Mac}98]{MacLane98}
	Saunders {Mac Lane}.
	\newblock {\em {Categories for the working mathematician. 2nd ed.}}
	\newblock New York, NY: Springer, 2nd ed edition, 1998.
	
	\bibitem[{May}68]{May1968}
	J.P. {May}.
	\newblock {Simplicial objects in algebraic topology.}
	\newblock {Princeton, N.J.-Toronto-London-Melbourne: D. van Nostrand Company,
		Inc. VI, 161 p. (1968).}, 1968.
	
	\bibitem[Mn{\"e}20]{MnevMin}
	N.~Mn{\"e}v.
	\newblock Minimal triangulations of circle bundles, circular permutations, and
	the binary {Chern} cocycle.
	\newblock {\em J. Math. Sci., New York}, 247(5):696--710, 2020.
	
	\bibitem[Mun00]{Munkres2000}
	James~R. Munkres.
	\newblock {\em Topology.}
	\newblock Upper Saddle River, NJ: Prentice Hall, 2nd ed. edition, 2000.
	
	\bibitem[MY91]{MorYo1991}
	Bernard Morin and Masaaki Yoshida.
	\newblock The {K{\"u}hnel} triangulation of the complex projective plane from
	the view point of complex crystallography. {I}.
	\newblock {\em Mem. Fac. Sci., Kyushu Univ., Ser. A}, 45(1):55--142, 1991.
	
	\bibitem[Pal61]{Palais1961}
	Richard~S. Palais.
	\newblock On the existence of slices for actions of non-compact {Lie} groups.
	\newblock {\em Ann. Math. (2)}, 73:295--323, 1961.
	
	\bibitem[{Ser}10]{Sergeraert2010}
	Francis {Sergeraert}.
	\newblock {Triangulations of complex projective spaces.}
	\newblock In {\em {Contribuciones cient\'\i ficas en honor de Mirian Andr\'es
			G\'omez.}}, pages 507--519. Logro\~no: Universidad de La Rioja, Servicio de
	Publicaciones, 2010.
	
	\bibitem[tDKP70]{DieckPuppe}
	Tammo tom Dieck, K.~H. Kamps, and Dieter Puppe.
	\newblock {\em Homotopietheorie}, volume 157 of {\em Lect. Notes Math.}
	\newblock Springer, Cham, 1970.
	
	\bibitem[Tsy83]{Tsygan1983}
	B.~L. Tsygan.
	\newblock The homology of matrix {Lie} algebras over rings and the {Hochschild}
	homology.
	\newblock {\em Russ. Math. Surv.}, 38(2):198--199, 1983.
	
\end{thebibliography}
\end{document}